\documentclass[amsmath,amssymb,hidelinks,onefignum,unsrt]{revtex4-2}
\usepackage{cleveref}
\usepackage{graphicx}

\newcommand{\pd}[2]{\frac{\partial #1}{\partial #2}}

\newcommand{\eps}{\epsilon}

\newcommand{\beq}{\begin{equation}}
\newcommand{\eeq}{\end{equation}}

\newcommand{\re}{\mbox{Re}}
\newcommand{\im}{\mbox{Im}}

\begin{document}

\title{Multiple-scale analysis of the simplest large-delay differential equation
}
\author{Gregory Kozyreff}
\affiliation{Optique Nonlin\'eaire Th\'eorique, Universit\'e libre de Bruxelles (U.L.B.), CP 231, Belgium}
\begin{abstract}
A delayed term in a differential equation reflects the fact that information takes significant time to travel from one place to another within a process being studied. Despite de apparent similarity with  ordinary differential equations, delay-differential equations (DDE) are known to be fundamentally different and to require a dedicate mathematical apparatus for their analysis. Indeed, when the delay is large, it was found that they can sometimes be related to spatially extended dynamical systems. The purpose of this paper is to explain this fact in the simplest possible DDE by way of a multiple-scale analysis. We show the asymptotic correspondence of that linear DDE with the diffusion equation. This partial differential equations arises from a solvability condition that differs from the ones usually encountered in textbooks on asymptotics: In the limit of large delays, the leading-order problem is  a map and secular divergence at subsequent orders stem from forcing terms in that map.
\end{abstract}
\date{\today}
%\pacs{05.65.+b, 47.54.-r, 47.11.St, 89.75.Kd}
\maketitle

%------------------------------------------------
\section{Introduction}
%------------------------------------------------
%The introduction introduces the context and summarizes the
%manuscript. It is importantly to clearly state the contributions of
%this piece of work. The next two paragraphs are text filler,
%generated by the \texttt{lipsum} package.

In mathematical modelling, to be able to describe a physical, chemical, or an automated process by a lumped-element model, \textit{i.e.}, by a finite set of time-dependent variables, is a promising starting point for fruitful analysis. Occasionally, the interaction between elements of the model takes place with a time delay that cannot be neglected. In such a case, the appropriate mathematical problem to be solved typically consists of one ore more Delay Differential Equations (DDEs), also called functional differential equations, as opposed to ordinary ones. The consequences can be dramatic: 
first order DDEs can have periodic or even chaotic solutions~\cite{Mackey-1977,Ikeda-1979,Farmer-1982,Erneux-2004,Sprott-2007,Mueller-2018,Banerjee-2018}. Hence,
despite their similarity in writing, DDEs are fundamentally different from ODEs and require specific consideration. It is easy to observe, for instance, that linear, constant-coefficient DDEs  generally admit an infinity of independent exponential solutions. They may thus be regarded as ODEs of infinite order~\cite{BellmanCooke}. In 1996, Giacommelli and Politi went one step further by showing the asymptotic equivalence  of some DDEs to partial differential equations when the delay is large~\cite{Giacomelli-1996,Giacomelli-1998}. This result followed earlier numerical simulations~\cite{Arecchi-1992} pointing to the relevance of a 2D representation of the solutions of large-delay DDEs, see \Cref{fig:DDE}(a). This ``spatio-temporal equivalence'' has been confirmed in several subsequent works~\cite{Wolfrum-2006,Wolfrum-2010,Yanchuk-2017} and the purpose of this paper is to discuss it  in the simplest possible setting, namely the equation
%\beq
%t_i u'(\tilde t)+u(\tilde t)=r u(\tilde t-t_d),
%\eeq
%in the limit
%\beq
%T=t_d/t_i\gg1.
%\eeq
%
%In this paper, we wish to explain the above equivalence in the simplest possible equation, namely
\begin{align}
y'(t)+y(t)&=ry(t-T), &T&\gg1,%, & y(t<0)&=\phi(t).
\label{simple}
\end{align} 
with  initial data $y(t)=\psi(t)$ in the range $t\in [-T,0]$.
The interpretation of this equation is straightforward: it models a linear system whose internal dynamics is governed by the left hand side, with a simple exponential decay, and which is subjected to a delayed feedback with strength $r$. Here, $T$ is the delay normalised by the eigenvalue %decay rate
 of the isolated system.  \Cref{simple} appears in recent developments as the deterministic part of a Langevin equation
\beq
y'(t)+y(t)=ry(t-T)+\mu\xi(t),
\label{noisy}
\eeq
where $\xi(t)$ is a noisy forcing term. The latter can model Brownian motion with a memory effect~\cite{Budini-2004,Munakata-2009} or delayed control~\cite{Ando-2017}. More recently, it was proposed to study chaotic diffusion mediated by a nonlinear DDE~\cite{Albers-2022}.
 \Cref{simple} is also a special case of
\beq
y'(t)=a y(t) + b y(t-T)
\eeq
which has been studied in detail as one of the simplest DDE~\cite{Hayes-1950,Stepan-1989,Michiels-2007}.

\vspace{.2cm}
Concretely, we will derive from \cref{simple}  the multiple-scale asymptotic approximation
\beq
y\sim r^s e^{-\ln(r)^2z/2}
 \times Y(s,z),
 \label{sol}
 \eeq
  with $s\sim t/(T+1)$, $z=t/T^3$ and
\begin{align}
\pd Yz&=\frac12\pd{^2Y}{s^2}, &Y(s,z)&=Y(s-1,z).
\label{intro:diffusion}
\end{align}

This calculation serves three pedagogical purposes:
\vspace{.2cm}

\begin{itemize}
\item to highlight the general interest of investigating the large-$T$ limit in DDEs,
\item to detail the inner workings of the method of multiple scales in that particular framework,
\item to provide instructors with a new, easily workable application of that method, beyond the usual perturbed harmonic oscillator, in the  spirit of~\cite{Edwards-2000}. 
\end{itemize}

\vspace{.2cm}

%--------------------------------------------------
\subsection{Application of the multiple-scale approximation}
\begin{figure}
\centering
\includegraphics[width=\textwidth]{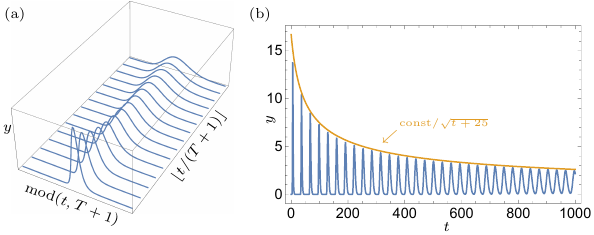}
%\caption{(a) Numerical solution of \cref{simple}  with $r=1$, $T=30$, and  initial condition $y(t)=20 e^{-(t+25)^2}$ in the range $-T<t<0$. The orange curve is the envelope of the peaks predicted by \cref{maxpeak}. (b) Pseudo spatio-temporal plot of the numerical solution, demonstrating the $T+1$ period of recurrence and the diffusive behaviour. }\label{fig:DDE}
\caption{Numerical solution of \cref{simple}  with $r=1$, $T=30$, and  initial condition $y(t)=20 e^{-(t+25)^2}$ in the range $-T<t<0$. (a) Pseudo spatio-temporal plot, demonstrating the $T+1$ period of recurrence and the diffusive behavior. (b) Usual representation as a function of $t$.}\label{fig:DDE}
\end{figure}

\begin{figure}
\centering
\includegraphics[width=8cm]{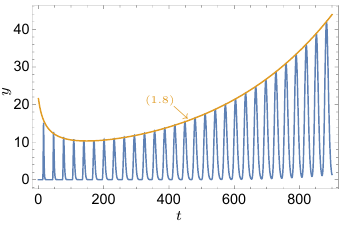}
\caption{Numerical solution of \cref{simple}  with $r=1.1$, $T=30$, and  initial condition $y(t)=20 e^{-(t+T/2)^2}$ in the range $-T<t<0$. The orange curve is the envelope of the peaks predicted by \cref{maxpeakr}. }\label{fig:DDEr}
\end{figure}

%--------------------------------------------------
Before deriving \cref{sol,intro:diffusion}, let us demonstrate their usefulness as a tool to analyze \cref{simple}. Recall first that an initial condition $\delta(s+s_0)$ at $z=z_0$ evolves under  the above diffusion equation as~\cite{Strauss-2007} 
\beq
Y(s,z)=\frac{e^{-(s+s_0)^2/2(z-z_0)}}{\sqrt{2\pi (z-z_0)}}.
\label{gaussian}
\eeq
From this, we may deduce,  when $r=1$, that a peaked initial condition will recur periodically with a flattening profile and a peak value that decreases in time as $1/\sqrt{t+c}$, where $c$ is an appropriate constant. Indeed, a gaussian initial condition 
\begin{equation*}
e^{-k(t+t_0)^2}
\end{equation*}
%\begin{align*}
%e^{-k(t+t_0)^2} &&\text{in \cref{simple}}
%\end{align*}
in \cref{simple} translates into $e^{-k(T+1)^2(s+s_0)^2}$ in \cref{intro:diffusion}. %, \textit{i.e.}, by periodicity in $s$, into $e^{-k(T+1)^2 s^2}$.
  It may thus be  assimilated, up to an appropriate factor, to a delta function in the $s$ variable in the limit $T\to\infty$. By the same token, any similarly localized initial condition can be asymptotically regarded as  proportional to a delta function on the $s$-scale. Therefore,  \cref{gaussian} applies. Furthermore, $t=-t_0$ corresponds to $z_0=-t_0/T^3$. Hence, evaluating \cref{gaussian} at $s=-s_0$, we deduce the following envelope for the maxima of $y(t)$:
\beq
%\max_{n(T+1)<t<(n+1)(T+1)} y(t)\propto
 \frac{\text{const}}{\sqrt{z-z_0}}\propto \frac{1}{\sqrt{t+t_0}}.
\label{maxpeak}
\eeq
 When $r\neq1$ the exponential factor in \cref{sol} is applied. One thus obtains the envelope
\beq
\text{const} \times \frac{r^{t\left(1/(T+1)+\ln(r) /T^3\right)}}{\sqrt{t+t_0}}.
\label{maxpeakr}
\eeq
We compare the above formulas with numerical solutions of \cref{simple} with $r=1$ and $r=1.1$ respectively  in \cref{fig:DDE,fig:DDEr}. \cref{fig:DDE}(a) confirms the relevance of the ``spatio-temporal'' representation of the solution, in which $\mod(t,T+1)\propto s$, $\lfloor x \rfloor$ denotes the integer part of $x$ and $\lfloor t/(T+1)\rfloor$ is a discrete variable on which the evolution is so slow that it is asymptotically equivalent to $z$. Next, \cref{fig:DDE}(b) and \cref{fig:DDEr} show the quantitative agreement of \cref{maxpeak,maxpeakr} with the numerical simulations. Especially striking is the non-monotonous behavior seen in the latter, which conforms to the intuition brought by \cref{intro:diffusion}. Indeed, diffusion promotes an initial flattening and attenuation of the peaks, before the amplification factor $r>1$ takes over.

Note that \cref{sol}  and \cref{gaussian}  actually make up an asymptotic approximation of the Green function of \cref{simple}, in an alternative way to the exact solution given in~\cite{Budini-2004}. This Green function is used to study the noisy extension \cref{noisy} of \cref{simple}.

% The outline is not required, but we show an example here.
The rest of the paper is devoted to the derivation of \cref{sol,intro:diffusion}. It starts  in \cref{sec:lin} with a standard linear analysis of \cref{simple}. The characteristic equation is derived and then simplified in the large-$T$ limit. This gives us insight on the relevant time scales for the multiple scale analysis, which we develop in \cref{sec:multiscale}. We carry out the calculation twice, as the method can be implemented in two slightly different but equally instructive  ways. 
 For the sake of simplicity, we initially focus on $r=1$ and
%We first solve the case $r=1$, as it is lightest from an algebraic point of view, and
 deal with general $r$ as a later step. Finally, we present our conclusions in  \cref{sec:discu}.
%\cref{sec:main}, our new algorithm is in \cref{sec:alg}, experimental
%results are in \cref{sec:experiments}, and the conclusions follow in
%\cref{sec:conclusions}.

%------------------------------------------------
\section{Linear spectrum}\label{sec:lin}
%------------------------------------------------

Let us look for a solution of (\ref{simple}) in the form of $y(t)=\exp\lambda t$. Assuming $r=1$ for simplicity, this yields the characteristic equation
\beq
\lambda+1=e^{-\lambda T}.
\label{eq:char}
\eeq 
Rearranging terms, we have
\begin{align}
\lambda&=-1+\frac uT, &
u e^{u}&=Te^T.
\end{align}
The equation $ye^y=x$ possesses a discrete infinity of complex solutions, denoted $W_n(x)$ and called ``ProductLog[n,x]'' in Mathematica. $W_n$ is the $n^\text{th}$ branch of the complex Lambert function~\cite{LambertWiki}. Hence, the spectrum of \cref{simple} is, exactly,
\beq
\lambda_n=-1+\frac{W_n\left(Te^T\right)}{T}.
\label{Lambert}
\eeq
%%%%%%%%%%%%%%%%%%%%%%%%%%%%
\begin{figure}
\centering
\includegraphics[width=8cm]{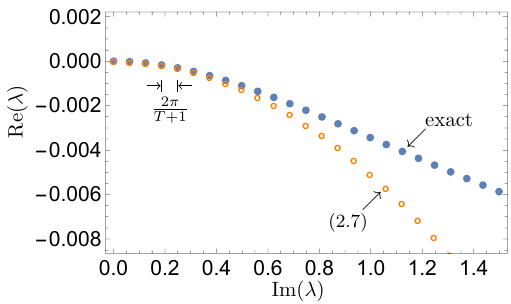}
\caption{The spectrum of \cref{simple} for $T=100$ (blue) compared to the approximation \cref{spectrumr=1} (orange circles). The separation between roots along the imaginary axis asymptotes to $2\pi/(T+1)$.}\label{figDDEspectrum}
\end{figure}
%%%%%%%%%%%%%%%%%%%%%%%%%%%%%
 \Cref{Lambert} is useful to draw the spectrum with a symbolic software but not very enlightening to anyone who is not an expert of the Lambert function. Fortunately, we can make significant progress in the limit $T\to\infty$. Writing $\lambda T=\sigma+i\omega$, \cref{eq:char} yields
 \begin{align}
\frac{ \sigma}{T}+1-e^{-\sigma }\cos\left(\omega \right)&=0,
&\frac{ \omega}{T}+e^{-\sigma }\sin\left(\omega \right)&=0. 
\label{eq:realim}
 \end{align}
This suggests the expansions $\sigma\sim \sigma_0+T^{-1}\sigma_1+T^{-2}\sigma_2+\dotsb$ and $\omega\sim \omega_0+T^{-1}\omega_1+\dotsb$. At leading order, we get
 \begin{align}
1-e^{-\sigma_0 }\cos\left(\omega_0 \right)&=0,
&e^{-\sigma_0 }\sin\left(\omega_0 \right)&=0,
 \end{align}
 which implies $\sigma_0=0$ and $\omega_0=2n\pi$, where $n$ is an integer. It is important to remark that $n$ must be assumed to be $O(1)$, so that $\omega/T$ can be treated as $O(1/T)$ in \cref{eq:realim}. To carry out the calculation to higher orders presents no difficulty and can be proposed as an exercise. One obtains
 \begin{align}
\sigma_1&=0, &\omega_1&=-\omega_0, &\sigma_2&=\frac{-\omega_0^2}{2},&\omega_2&=\omega_0, &\dotsb
 \end{align}
 Eventually,
 \begin{equation}
 \lambda\sim\frac1T\left[-\frac{2n^2\pi^2}{T^2}+2in\pi\left(1-\frac1T+\frac1{T^2}\right)+\dotsb\right]%\\
 \sim\frac1T\left(-\frac{2n^2\pi^2}{T^2}+\frac{2in\pi}{1+1/T}\right)=
 -\frac{2n^2\pi^2}{T^3}+\frac{2in\pi}{T+1}.
 \label{spectrumr=1}
 \end{equation} 
 The exact and approximate spectra are depicted in \cref{figDDEspectrum} and are found to be in good agreement  as long as the imaginary part is small, \textit{i.e} for mode numbers satisfying $2\pi n\ll T$.

On the other hand,  note that the diffusion equation \cref{intro:diffusion} has exponential solutions $\exp(\kappa z+i \Lambda s)$  provided that
\beq
\kappa = -\Lambda^2/2,
\eeq
while the periodic boundary condition in $s$ imposes $\Lambda=2n\pi$, with integer $n$. Using the definitions of $s$ and $z$, one thus obtains
$\exp\{[-2n^2\pi^2 /T^3+2in\pi /(T+1)]t\}$, which is consistent  with \cref{spectrumr=1}.

 The imaginary part of $\lambda$ points to an oscillatory evolution with approximate period $T+1$. Notice  that time in \cref{simple} is rescaled to the  intrinsic time scale of the isolated ($r=0$) dynamical system. With a more general time unit, this intrinsic time scale would  numerically  be given by $\bar t_i$, the delay by  $\bar T$, and the period of damped oscillations by $\bar T+\bar t_i$. In view of this,  the periodicity is fixed by the time required for information to be fed back into the dynamical system plus the time required to internally process it. Such a conclusion would be difficult to draw from the contemplation of \cref{Lambert} alone.

The analysis for $r\neq1$ follows the same pattern. We leave it as an exercise to show that, in the general case,
\begin{equation}
 \lambda\sim\frac1T\left[-\frac{\ln r}{2T^2}-\frac{2n^2\pi^2}{T^2}+\left(\ln r+2in\pi\right)\left(1-\frac1T+\frac{1+\ln r}{T^2}\right)+\dotsb\right].
 \label{spectrumwithr}
 \end{equation}

%------------------------------------------------
\section{Multiple-scale analysis}\label{sec:multiscale}
%------------------------------------------------
The linear stability analysis with $r=1$ reveals that the time scales for oscillations and damping become infinitely separated as $T\to\infty$, which suggests a multiple-scale analysis. In this section, we propose two slightly distinct implementations of the method. The first one is a two-time expansion of the solution, 
\beq
y(t)\sim Y(s,z),
\label{twotime}
\eeq
in which the fast time is a strained coordinate:
\beq
s=\left(1+\frac{a_1}{T}+\frac{a_2}{T^2}+\dotsb\right)\frac tT,
\label{eq:strain}
\eeq
and $z=t/T^3$. In that implementation, we  construct a solution whose  periodicity is strictly $1$ in $s$. The constants $a_1$ and $a_2$ are determined in the course of resolution, but they can in fact be guessed from the results of the linear stability analysis. 

The alternative is to pose a three-time ansatz
\beq
y(t)\sim Y(\tau,\eta,z),
\label{threetime}
\eeq
where $\tau=1/T$ and $\eta=t/T^2$ with no preconception of how  $f$ must depend on each time scale. Both approaches have their advantages and disadvantages. To introduce multiple time scales and convert a total derivative into a partial differential operator can be a deterring prospect to students who are exposed to the method for the  first time. This speaks in favor of minimizing the number of time scales and, hence, of choosing the ansatz (\ref{twotime}). On the other hand, the strained $s$ coordinate may appear pulled out of a hat. This can make (\ref{threetime})  more easily acceptable to those who prefer a step-by-step approach. Eventually, given the lightness of the calculations, it may be useful for the instructor to present the two variants of the calculation and, in so doing,  to demonstrate the robustness of the method. In order to keep the calculation down to its essential details, we focus first on the case $r=1$. We treat the case $r\neq1$ as a later step.

%-------------
\subsection{Two-time calculation ($r=1$)}
Assuming that $y$ is asymptotically given by the ansatz~(\ref{twotime}), differentiation with respect to $t$ yields, by the chain rule,
\beq
y'(t)\sim \left(1+\frac{a_1}{T}+\frac{a_2}{T^2}+\dotsb\right)\frac 1T\pd{Y}{s}+\frac1{T^3}\pd{Y}{z}.
\eeq
On the other hand, the delayed term becomes
\begin{multline}
y(t-T)\sim Y\left(s-1-\frac{a_1}{T}-\frac{a_2}{T^2}-\dotsb,z-\frac1{T^2}\right)\\
\sim Y\left(s-1,z\right)-\frac{a_1}{T}\pd{}{s}Y\left(s-1,z\right)+\frac1{T^2}
\left[\frac{a_1^2}{2}\pd{^2}{s^2}-a_2\pd{}{s}-\pd{}{z}\right]Y\left(s-1,z\right)+\dotsb
\end{multline}
Hence, \cref{simple} is transformed into
\begin{multline}
\left(\frac1T+\frac{a_1}{T^2}\right)\pd{Y(s,z)}{s}+Y(s,z)
\sim  Y\left(s-1,z\right)-\frac{a_1}{T}\pd{}{s}Y\left(s-1,z\right)\\
+\frac1{T^2}
\left[\frac{a_1^2}{2}\pd{^2}{s^2}-a_2\pd{}{s}-\pd{}{z}\right]Y\left(s-1,z\right) +O\left(\frac1{T^3}\right).
\label{simple2time}
\end{multline}
Expanding $Y$ as $Y_0+T^{-1}Y_1+T^{-2}Y_2+\dotsb$, one obtains, at $O(1)$,
\beq
Y_0(s,z)=Y_0\left(s-1,z\right).
\label{2time:order0}
\eeq
%Posing $f_0(s,z)=r^sF_0(s,z)$, this yields
%\beq
%F_0(s,z)= F_0\left(s-1,z\right),
%\label{2time:order0bis}
%\eeq
%so that $F_0$ is periodic with period $1$ in $s$. 
so that $Y_0$ is periodic with period $1$ in $s$. 
%
%Hence, up to a $r^s$ factor (which disappears if $r=1$), the solution is periodic of period $1$ in $s$.
%
%The leading order problem thus asserts that, up to the $r^s$ factor (which disappears if $r=1$), the solution is periodic of period $1$ in $s$.
% At this stage, no other constraint exists on $Y_0$: any combination of Fourier component  $\exp(2in\pi s)$ is allowed.
  Next, collecting $O(1/T)$ terms and using \cref{2time:order0}, we find
\beq
Y_1(s,z)= Y_1\left(s-1,z\right)-\left(1+a_1\right) \pd{}{s}Y_0(s,z) .
\eeq
The forcing term, being periodic, causes a secular divergence of $f_1$:
\beq
Y_1(s+j,z)= Y_1\left(s,z\right)-j\left(1+a_1\right)  \pd{}{s}Y_0(s,z).
\label{discretegrowth}
\eeq
Hence, the asymptotic ordering $T^{-1}Y_1\ll Y_0$ breaks down for $j=O(T)$, unless we set
\beq
a_1=-1.
\eeq
At $O(T^{-2})$, we then find
\begin{equation}
Y_2(s,z)= Y_2\left(s-1,z\right)
+\left(1-a_2\right)\pd{Y_0}{s}+\frac12\pd{^2Y_0}{s^2}-\pd{Y_0}{z}
,
\end{equation}
where $Y_0$ is evaluated at $(s,z)$. 
Following the same reasoning as at the previous order, the solvability condition is
\beq
\pd{Y_0}{z}=\left(1-a_2\right)\pd{Y_0}{s}+\frac12\pd{^2Y_0}{s^2},
\label{nearlythere}
\eeq
At this stage, it may seem that the constant $a_2$ is still free. However, if we insist that the periodic component of the solution is strictly $1$ in $s$, then $Y_0$ is a combination of functions of the form $\exp(\kappa z +2in\pi s)$, where, according to \cref{nearlythere}, $\kappa= 2in\pi (1-a_2)-2n^2\pi^2$. Now, the imaginary part of $\kappa$ perturbs the period of the solution and to avoid this, we make it vanish:
\beq
a_2=1,
\eeq
which agrees with \cref{spectrumr=1}. We thus obtained \cref{intro:diffusion}.

%-------------
\subsection{Three-time calculation}\label{sec:3time}
If we assume the multi-time ansatz  \cref{threetime}, then
\beq
y'(t)\sim\frac1T \pd Y\tau+\frac1{T^2}\pd{Y}{\eta}+\frac{1}{T^3}\pd{Y}{z}
\eeq
and
\begin{multline}
y(t-T)\sim Y\left(\tau-1,\eta-T^{-1},z-T^{-2}\right)
\\
\sim
Y\left(\tau-1,\eta,z\right)
%\\
-\frac1T \pd{}{\eta} Y\left(\tau-1,\eta,z\right)-\frac1{T^2}\left(\pd{}{z} -\frac12\pd{^2}{\eta^2}\right)Y\left(\tau-1,\eta,z\right).
\label{3-timedelay}
\end{multline}
We thus get to solve
\begin{equation}
%\frac1T \pd Y\tau+\frac1{T^2}\pd{Y}{\eta}+\frac{1}{T^3}\pd{Y}{z}+Y=\\
\left(1+\frac1T \pd{}\tau+\frac1{T^2}\pd{}{\eta}+\frac{1}{T^3}\pd{}{z}\right)Y\left(\tau,\eta,z\right)=
\left(1-\frac1T\pd{}{\eta}-\frac1{T^2}\pd{}{z} +\frac1{2T^2}\pd{}{\eta^2}+\dotsb\right) Y\left(\tau-1,\eta,z\right).
\end{equation}
Expanding $Y$ again as $Y_0+T^{-1}Y_1+\dotsb$, we obtain
\beq
Y_0(\tau,\eta,z)=Y_0(\tau-1,\eta,z),
\label{3time:order0}
\eeq
\textit{i.e.},  the same map as before, with $s$ replaced by $\tau$. At $O(T^{-1})$, we have
\begin{equation}
Y_1(\tau,\eta,z)-Y_1(\tau-1,\eta,z)=-\pd{}{\tau}Y_0(\tau,\eta,z)-\pd{}{\eta}Y_0(\tau-1,\eta,z)
=-\left(\pd{}{\tau}+\pd{}{\eta}\right)Y_0(\tau,\eta,z).
\label{eq4Y1}
\end{equation}
Since the right hand side is a period-1 function in $\tau$, solvability requires that it vanishes:
\beq
\pd{Y_0}{\tau}+\pd{Y_0}{\eta}=0.
\eeq
The general solution is, simply,
\beq
Y_0(\tau,\eta,z)=Y_0(\tau-\eta,z),
\eeq
or, equivalenty,
\begin{align}
Y_0(\tau,\eta,z)&=Y_0(S,z), &S=\tau-\eta=\frac tT\left(1-\frac1T\right).
\end{align}
The equation for $Y_1$ then implies that it is periodic in $\tau$. Next, at $O(T^{-2})$, and taking into account the periodicity of $Y_0$ and $Y_1$, the problem for $Y_2$ is
\begin{equation}
Y_2(\tau,\eta,z)-Y_2(\tau-1,\eta,z)=-\left(\pd{}{\tau}+\pd{}{\eta}\right)Y_1(\tau,\eta,z)+\left(\pd{}{S}
+\frac1{2}\pd{}{S^2}-\pd{}{z} \right)Y_0(S,z).
\end{equation}
Secular divergence of $Y_2$ is avoided by making the right hand side vanish:
\beq
\pd{Y_1}{\tau}+\pd{Y_1}{\eta}=\pd{Y_0}{S}
+\frac1{2}\pd{^2Y_0}{S^2}-\pd{Y_0}{z} .
\label{SC:primary}
\eeq
This is an equation for $Y_1$, with general solution
\beq
Y_1(\tau,\eta,z)=Y_1(S,z)+\eta\left(\pd{Y_0}{S}  +\frac1{2}\pd{^2Y_0}{S^2}-\pd{Y_0}{z}\right).
\eeq
There is therefore a secular divergence in $\eta$ of $Y_1$ and to avoid this, we have a solvability condition on the solvability condition:
\beq
\pd{Y_0}{z}-\pd{Y_0}{S} =\frac1{2}\pd{^2Y_0}{S^2}.
\label{SC:secondary}
\eeq
Finally, if we write 
\begin{align}
Y_0(S,z)&=\phi(s,z), &s=S+z=\tau-\eta+z=\frac tT\left(1-\frac1T+\frac1{T^2}\right),
\end{align}
%where $s$ is the same rescaled time variable as in the previous section,
 we obtain
\beq
\pd{\phi}{z} =\frac1{2}\pd{^2\phi}{s^2},
\label{threetime:diff}
\eeq
\textit{i.e.} the desired result. Note that the fact that $Y_0(\tau,\eta,z)$ is of period $1$ in $\tau$ means that $\phi$ is of period $1$ in $s$, in full consistency with what precedes.

\vspace{.1cm}

\textbf{Remark} One may be tempted to  purely and simply set $Y_1=0$, after noting (i) that the initial problem is linear and (ii) that \cref{eq4Y1} is homogenous after applying the solvability condition. Indeed, this would expedite the derivation of \cref{threetime:diff}. However, it would conceal the structure of nested solvability conditions appearing at $O(T^{-2})$. While there is nothing wrong in seeking the shortest route to the answer, especially in the eyes of an applied mathematician, one should bear in mind that $Y_1$ should, in all generality, be retained. This correction may be required to properly handle $O(T^{-1})$ terms in the initial condition or weak nonlinearities. Interestingly, while there is only a single equation to solve at  $O(T^{-2})$, more than one solvability conditions can be extracted from it: the main one, \cref{SC:primary}, and the secondary one, \cref{SC:secondary}.

%%//////////////////////////////////
%\newpage
%\appendix
\subsection{The case $r\neq1$}
We now revise the calculation in the more general case $r\neq1$. Here again, the two implementations yield the same result. We limit ourselves to the ``two-time'' calculation and assume
\beq
y\sim f(s,z),
\eeq
where $s$ is given by \cref{eq:strain}. This time,  \cref{simple} is transformed into
\begin{multline}
\left(\frac1T+\frac{a_1}{T^2}\right)\pd{f(s,z)}{s}+f(s,z)
\sim r\bigg\{f\left(s-1,z\right)-\frac{a_1}{T}\pd{}{s}f\left(s-1,z\right)\\
+\frac1{T^2}
\left[\frac{a_1^2}{2}\pd{^2}{s^2}-a_2\pd{}{s}-\pd{}{z}\right]f\left(s-1,z\right)\bigg\}+O\left(\frac1{T^3}\right).
\end{multline}
Expanding $f$ as $f_0+T^{-1}f_1+T^{-2}f_2+\dotsb$, the  $O(1)$ problem is now
\beq
f_0(s,z)= r f_0\left(s-1,z\right).
\label{rneq1:order0}
\eeq
Posing $f_0(s,z)=r^sF_0(s,z)$, this yields
\beq
F_0(s,z)= F_0\left(s-1,z\right),
\label{rneq1:order0bis}
\eeq
so that $F_0$ is periodic with period $1$ in $s$. Next, collecting $O(1/T)$ terms and using \cref{rneq1:order0,rneq1:order0bis}, we find
\beq
f_1(s,z)= r f_1\left(s-1,z\right)-\left(1+a_1\right)r^s\left[\ln(r) F_0(s,z)+\pd{}{s}F_0(s,z)\right].
\eeq
The terms between bracket, being periodic, cause a secular divergence of $f_1$. Indeed, letting $f_1(s,z)=r^sF_1(s,z)$, one rapidly finds that
\beq
F_1(s+j,z)= F_1\left(s,z\right)-\left(1+a_1\right) j \left[\ln(r) F_0(s,z)+\pd{}{s}F_0(s,z)\right].
\eeq
Hence, irrespective of the factor $r^s$, $T^{-1}f_1$ ceases to be small compared to $f_0$ when $j=O(T)$. To prevent this, we thus set
\beq
a_1=-1.
\eeq
At $O(T^{-2})$, we have
\begin{equation}
f_2(s,z)= r f_2\left(s-1,z\right)
+r^s\Bigg\{
\left[\frac{\ln(r)^2}2+\ln(r)\left(1-a_2\right)\right]F_0(s,z)
+\left(1+\ln(r)-a_2\right)\pd{F_0}{s}+\frac12\pd{^2F_0}{s^2}-\pd{F_0}{z}
\Bigg\}.
\end{equation}
The solvability condition is now
\beq
\pd{F_0}{z}=l_0 F_0(s,z)
+l_1\pd{F_0}{s}+\frac12\pd{^2F_0}{s^2},
\eeq
where $l_0=\frac{\ln(r)^2}2+\ln(r)\left(1-a_2\right)$ and $l_1=1+\ln(r)-a_2$. The term multiplied by $l_1$ causes a change in the periodicity of the solution and is therefore set to zero:
\beq
a_2=1+\ln(r).
\eeq
This agrees with \cref{spectrumwithr}. Hence, $l_0=\frac{-\ln(r)^2}2$. It only remains to make a small change of variable, namely to write $F_0=e^{l_0 z}Y(s,z)$ to obtain \cref{{intro:diffusion}}.

%---------------------------------------------
\section{Discussion}\label{sec:discu}

The derivation presented in this paper hints to the great generality of the diffusion equation \cref{intro:diffusion} as the linear backbone of large-delay differential equations. Not all DDE develop the multiple-scale structure presented here (see below) but when they do, diffusion is to be expected from the Taylor expansion of the delayed term with the appropriate scales: see  the  $O(T^{-2})$ terms in the right hand side of \cref{3-timedelay}.

\vspace{.2cm}
An interpretation of the delay in a DDE like \cref{simple} is that the output of a given system undergoes some kind of propagation before being fed back. The form of the delayed term, $r y(t-T)$, is strongly suggestive of a hyperbolic PDE as a mediator of this feedback, a point of view that was emphasized in~\cite{Buenner-1997}. 
%Hence, it is somewhat surprising, physically, that the parabolic PDE \cref{intro:diffusion} should arise out of a hyperbolic one. That it is so results from 
Why, then, should a parabolic PDE such as \cref{intro:diffusion}  arise out of a hyperbolic one? The answer to this question lies in 
the fact that the system subjected to feedback is itself dissipative, being described in isolation by the left hand side of \cref{simple}. Loosely speaking, a complete feedback cycle %thus
 includes a dispersion-less propagation of duration $T$ followed by an attenuation of unit duration. In the spectral domain, the left hand side of \cref{simple} is a low-pass filter. This is where information is degraded, as happens in diffusive processes. In this regard, one should bear in mind that, numerically, $T$ is the ratio of the delay to the internal dissipation time of the system. As we wrote earlier, in a general unit system, the delay is numerically given by $\bar T$ and the internal dynamics is characterized by $\bar t_i$, with $T=\bar T/\bar t_i$. It follows that the parameter $T$ does not appear as a result of the delay alone. It exists thanks to both the delay \emph{and} the internal time scale. Even if the latter is short compared to the former, the dissipation process occurring on the internal timescale is not to be neglected. 
Another facet of this question is given by the picture of the linear spectrum in \cref{figDDEspectrum}. Here,  the plot of the real part of the eigenvalues as a function of their imaginary part is analogous to  dispersion relations found in pattern formation,  where the exponential growth rate of a perturbation is plotted as a function of its wave number~\cite{Cross-1993,Hoyle-2006,Pismen-2006}. In that context, the existence of a small band of wave numbers with near-zero, maximum growth rate generically leads to diffusive terms in amplitude equations. Presently, the spectrum is discrete but tends to a continuum as $T\to\infty$. Hence, with a trained eye, diffusion can directly be anticipated at the sight of \cref{figDDEspectrum}.

\vspace{.2cm}
A  general technical feature exemplified by the present calculation is that the leading order problem is a map, see \cref{2time:order0,3time:order0,rneq1:order0}. This is because the main time scale is asymptotically set by the delay, and time derivatives of functions evolving on this time scale  are $O(T^{-1})$. In problems where the period is not asymptotically given  by the delay,  a map is \emph{not} expected as the leading order problem, even if the delay is large. Consider for example  Minorsky's equation~\cite{Pinney-1958}
\begin{align}
y''(t)+\eps y'(t)+\Omega^2 y(t)&=-by'(t-T)+\eps c y'(t-T)^3, &\eps&\ll1, &T=O\left(\eps^{-2}\right),
\end{align} 
which displays a  Hopf bifurcations with an $O(1)$ frequency. Here $\Omega$ sets the oscillation period and the leading-order problem of the multiple-scale analysis is the familiar harmonic oscillator. The linear spectrum is again densified by the largeness of $T$ but it now displays a maximum with near vanishing $\re(\lambda)$ in the vicinity of $\im(\lambda)=\Omega$. This is akin to a Turing bifurcation in spatially extended dynamical system. Following the preceding discussion, diffusion is therefore expected in the amplitude equation. Taking nonlinear terms into account one eventually obtains a Ginzburg-Landau equation (see~\cite{Erneux-2009} for a detailed derivation.) In any case, whether the leading order problem is a map or a harmonic oscillator, the general procedural idea of the multiple scale analysis is the same: identify and kill secular divergences. %Here, the  divergence is discrete, see \cref{discretegrowth}. %This, together with the remark at the end of \cref{sec:3time} may be worth underlining when teaching the method of multiple sclae. 

\begin{figure}
\centering
\includegraphics[width=8cm]{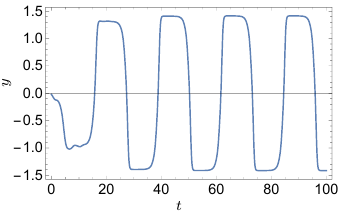}
\caption{Numerical solution of \cref{counterexample} with $T=10$ and $y(t)=0.1\sin(t)-0.02\cos(3t) $ for $t<0$ as an example of a non-diffusive long-term behavior.}\label{fig:square}
\end{figure}

The present paper may give the reader  the false impression  that  multiple-scale analysis is the method of choice to study all DDE in the large-delay limit. This is not the case and the following equation is a simple counter-example:
\beq
-y'(t)+y(t)=y(t-T)+y^3(t).
\label{counterexample}
\eeq
Compared to \cref{simple}, we have simply changed the sign of the time derivate and added a nonlinear term to avoid blow up. A numerical simulation is shown on \cref{fig:square}. Independently of the initial condition, the solution of this equation asymptotes to  sustained square wave oscillations of period $2T$ that display abrupt switching between approximately $\sqrt 2$ and $-\sqrt2$. This dynamical behavior is obviously not compatible with \cref{intro:diffusion}.  To study such a limit cycle, the method of matched asymptotic expansions appears more appropriate~\cite{Fowler-1982,Nizette-2004}.

\acknowledgments
G.K. is supported by the Fonds de la Recherche Scientifique - FNRS (Belgium.) 

%\bibliographystyle{siamplain}
%\bibliography{../siamart_220329/DDEbib}

%apsrev4-2.bst 2019-01-14 (MD) hand-edited version of apsrev4-1.bst
%Control: key (0)
%Control: author (8) initials jnrlst
%Control: editor formatted (1) identically to author
%Control: production of article title (0) allowed
%Control: page (0) single
%Control: year (1) truncated
%Control: production of eprint (0) enabled
%

\end{document}